\documentclass[a4paper,3p,11pt,sort&compress]{elsarticle}
\journal{Comput Methods Appl Mech Eng, final version: 10.1016/j.cma.2023.116488}
\date{}

\RequirePackage{snapshot}

\usepackage[Symbol]{upgreek}
\usepackage{mathptmx}
\usepackage[T1]{fontenc}
\usepackage[latin9]{inputenc}
\usepackage{amsmath}
\usepackage{amssymb}
\usepackage{graphicx}
\usepackage{setspace}
\usepackage{esint}
\PassOptionsToPackage{normalem}{ulem}
\usepackage{ulem}
\usepackage{marvosym}
\usepackage{multirow}
\usepackage{soul} 
\usepackage{physics}
\usepackage{xfrac}
\usepackage{dsfont}
\usepackage[hidelinks]{hyperref}
\usepackage{nomencl}
\usepackage{adjustbox}
\usepackage{mathtools}
\usepackage[hang,flushmargin]{footmisc}
\usepackage{subfig}
\usepackage{bm}
\usepackage{color}
\usepackage{epstopdf}

\soulregister\cite7
\soulregister\ref7
\soulregister\pageref7

\def\*#1{\mathbf{#1}}

\newdefinition{rmk}{Remark}

\makeatletter

\newcommand{\m}[1]{\mathbf{#1}}                     % matrix
\newcommand{\vt}[1]{\mathbf{#1}}                    % vector
\newcommand{\vg}[1]{\bm{#1}}                        % greek vector
\newcommand{\T}[0]{\mathrm{T}}                      % transpose
\newcommand{\cdotn}[0]{\!\cdot\!}                   % cdot with smaller spaces
\newcommand{\uk}{u}
\newcommand{\ukh}{u_h}
\newcommand{\tf}{v}
\newcommand{\tfh}{v_h}

\newcommand{\til}[1]{\widetilde{#1}}                    % tilde
\newcommand{\mhat}[1]{\hat{#1}}                         % error symbol
\newcommand{\dotP}[0]{\scalebox{1}{\textbullet}}        % larger \cdot

\DeclareMathSymbol{\shortminus}{\mathbin}{AMSa}{"39}
\newcommand{\shortm}{\negthinspace\shortminus\negthinspace}

\DeclareMathAlphabet\mathbfcal{OMS}{cmsy}{b}{n}         % bold mathcal
\DeclareMathAlphabet\mathcal{OMS}{cmsy}{n}{n}           % bold mathcal
\graphicspath{{figures/}}

\makeatother
\PassOptionsToPackage{numbers,sort&compress}{natbib}

\makenomenclature
\begin{document}

\makeatletter
\def\bm@pmb@#1{{%
  \setbox\tw@\hbox{$\m@th\mkern.3mu$}%
  \mathchoice
    \bm@pmb@@\displaystyle\@empty{#1}%
    \bm@pmb@@\textstyle\@empty{#1}%
    \bm@pmb@@\scriptstyle\defaultscriptratio{#1}%
    \bm@pmb@@\scriptscriptstyle\defaultscriptscriptratio{#1}}}
\makeatother

\title{Stabilized finite elements for the solution of the\\ Reynolds equation considering cavitation}

\author[cimne]{Hauke~Gravenkamp\corref{cor1}}
\ead{hgravenkamp@cimne.upc.edu}

\author[ovg]{Simon~Pfeil}

\author[cimne,upc]{Ramon~Codina}

\address[cimne]{International Centre for Numerical Methods in Engineering (CIMNE), 08034 Barcelona, Spain}

\address[ovg]{Otto von Guericke University Magdeburg, Institute of Mechanics, 39106 Magdeburg, Germany}

\address[upc]{Universitat Polit\`ecnica de Catalunya, 08034 Barcelona, Spain}
\cortext[cor1]{Corresponding author }

\begin{abstract}
The Reynolds equation, combined with the Elrod algorithm for including the effect of cavitation, resembles a nonlinear convection-diffusion-reaction (CDR) equation. Its solution by finite elements is prone to oscillations in convection-dominated regions, which are present whenever cavitation occurs. We propose a stabilized finite-element method that is based on the variational multiscale method and exploits the concept of orthogonal subgrid scales. We demonstrate that this approach only requires one additional term in the weak form to obtain a stable method that converges optimally when performing mesh refinement.
\end{abstract}
\begin{keyword}
    hydrodynamic bearings, Reynolds equation, stabilized finite elements, variational multiscale method, orthogonal subgrid scales
\end{keyword}
\maketitle

\section{Introduction}\noindent
The Reynolds equation is essential in modeling hydrodynamic lubrication processes in bearings \cite{Reynolds1886}. The basic form of this partial differential equation (PDE) is derived from the Navier-Stokes equations, describing the special case of pressure generation in thin fluid films. 
Early approaches to the numerical solution of the Reynolds equation date back several decades and include conventional techniques, such as the finite-element method \cite{Reddi1969}, finite differences \cite{Gnanadoss1964}, and the finite-volume method \cite{Arghir2002}. A semi-analytical technique has recently been described for obtaining efficient solutions under some simplifying assumptions \cite{Pfeil2020a}.

It is well known that the Reynolds equation, in its simplest version, can lead to unphysical negative pressure results in regions where the physical pressure becomes small enough for cavitation to occur. As a consequence, different variants have been proposed to obtain a more realistic model that includes effects due to cavitation. A trivial solution consists in setting any negative values to zero in a postprocessing step and, as a consequence, violating mass conservation. This method is often referred to as G\"umbel or half-Sommerfeld condition, see, e.g., \cite{Li2017b}. More physically accurate cavitation models enforce mass conservation through either the boundary conditions or the differential equation, in both cases leading to a nonlinear boundary value problem (BVP). In the former case, the computational domain is limited to the non-cavitated region (pressure zone) or the two flow regimes are interpreted as two separate domains coupled through interface conditions. The nonlinearity arises from the fact that the locations of the pressure and cavitation zones are initially unknown. This concept leads to the Swift-Stieber (or Reynolds) boundary condition~\cite{Swift1932,Stieber1933} and -- in a more sophisticated variant -- to the Jakobsson-Floberg-Olsson (JFO) conditions~\cite{Jakobsson1957,Floberg1957,Olsson1965}. See, e.g.,~\cite{Schweizer2008,Schweizer2008a} for analyses and numerical implementations of these two approaches. The other class of cavitation models incorporates the entire fluid film, i.e., pressure as well as cavitation zones, into one common computational domain, usually in combination with linear boundary conditions. The Reynolds equation is then formulated as a globally valid nonlinear differential equation considering the different physical properties of both regimes. The Elrod (or Elrod-Adams) algorithm~\cite{Elrod1974,Elrod1981,Fesanghary2011}, which is the most prevalent cavitation model, and the bi-phase approach~\cite{Feng1986,Zeidan1989} fall into this category. These cavitation algorithms still satisfy (and are motivated by) the interface conditions stipulated by the JFO model without enforcing them explicitly.
The cavitation model used in the study at hand is based on the assumptions proposed by Kumar and Booker~\cite{Kumar1991}.
Our current work will focus on the numerical treatment of this existing model; hence, we will not discuss in great detail the peculiarities of the modeling approach. Section~\ref{sec:problem} will provide a brief summary of the BVP; for a more detailed outline of its derivation, the reader may consult \cite{Pfeil2023}, in addition to the original works cited above.
Note that this cavitation model, like many others, is usually classified as a variant of the Elrod algorithm, as the assumptions are similar (although not identical) to the original work by Elrod and Adams.

The numerical solution of the Reynolds equation is prone to instabilities and severe oscillations, which can severely hinder the convergence of the nonlinear terms and require infeasibly small mesh sizes to obtain realistic solutions.
A conventional approach to improving stability and convergence relies on including an additional unphysical diffusion term that is nonzero in the cavitation domain.
This method is referred to as artificial diffusion or artificial dissipation and is used in, e.g.,~\cite{Shi2002,vanostayen2009,Alakhramsing2015}. 
Alternatively, numerical diffusion can be introduced through upwind schemes.
In the context of finite element models, this is achieved through Petrov-Galerkin formulations~\cite{Hajjam2007,Lengiewicz2014} -- in particular the streamline upwind Petrov-Galerkin (SUPG) method \cite{Habchi2012,Liu2020} -- where the test functions are modified, or by a special quadrature approach~\cite{Bertocchi2013} where the positions of the integration points are shifted depending on the flow direction. 
In finite difference or finite volume schemes, upwind differences are used, meaning that the convective term is discretized with a backward difference relative to the flow direction instead of a central difference~\cite{Vijayaraghavan1989,Shyu2008,Ausas2009}.
The amount of diffusion needs to be large enough to stabilize the solution but also small enough to avoid overly diffusive results.
Since finer meshes require less diffusion to achieve stability, the tuning parameters in the artificial diffusion method are usually defined as functions of the element size. Similar mesh-dependent behavior is obtained when employing upwind schemes.
In this study, we propose a stabilization technique based on the variational multiscale (VMS) method that is able to achieve optimal error convergence under mesh refinement and, to our knowledge, has not been applied to the Reynolds equation as of yet.
We may stress that the objective of this study is not the development of a new cavitation model but rather the design of a stabilization technique for an existing one.

The oscillatory behavior, as observed in the solution of the Reynolds equation, is also well-known from the analysis of the (linear) convection-diffusion-reaction (CDR) equation or the Navier-Stokes equations. In particular, instabilities are observed in `convection-dominated flows,' i.e., where convective terms (involving first-order spatial derivatives) are significantly larger than diffusive terms (derivatives of second order)~\cite{Codina2000a}. We will demonstrate that the Reynolds equation in the particular version discussed here can be cast into the form of a nonlinear and inhomogeneous CDR equation such that methods developed for the linear CDR equation can be adapted to our application. 
In particular, we devise a stabilized finite element method to prevent global oscillations and improve the convergence of the nonlinear problem. The proposed approach is based on the concept of the VMS method \cite{Hughes1995a,Hughes1998, Codina2017a}. The general idea of this framework is based on splitting the unknown solution into a contribution that is in the finite element space and a remainder that cannot be represented by the chosen discretization. Within this family of approaches, existing variants differ in the assumed model for this remainder, often referred to as \textit{subscales}. An overview of several methods for solving the CDR equation is presented in \cite{Codina1998}. 
An important subclass can be summarized as residual-based methods \cite{Lins2010}, which include stabilization terms involving the strong residual of the governing equation to ensure the error introduced by the stabilization vanishes as the element size approaches zero. 
We will use here an extension of this idea that involves the projection of the residual onto a space orthogonal to the finite element space. Details on this approach (now coined \textit{orthogonal subgrid scales}, OSGS) for the linear CDR equation were presented in \cite{Codina2000}.

A major advantage of the OSGS approach lies in the fact that it allows neglecting terms in the residual that are not essential for achieving stability while retaining optimal convergence of the method. This idea is sometimes referred to as \textit{term-by-term stabilization} \cite{Codina2008, Coppola-Owen2011, Castillo2017, Castillo2019}. Based on this concept, we are able to propose a stable and optimally converging method by including only one additional term in the discretized weak form. In addition to deriving the stabilized finite-element method, we will briefly discuss the application of a shock-capturing method that is beneficial in problems involving large gradients of the solution. 

In the following section, we will briefly summarize the BVP that will be the focus of this paper before discussing a suitable variational form in Section~\ref{variational}. The main idea -- the proposed stabilized finite element method for the Reynolds equation -- is explained in Section~\ref{stabilization}. In the ensuing, we discuss some essential numerical aspects, namely the computation of the projections (Section~\ref{projection}), the linearization of the discretized weak form (Section~\ref{linearization}), and the application of a shock-capturing algorithm (Section~\ref{shockCapturing}). Finally, Section~\ref{numex} presents three numerical studies demonstrating stable and optimally converging results.

\section{Problem statement}\label{sec:problem}\noindent
The model assumed here has been described in much detail in \cite{Pfeil2023} and the references therein. Here, it shall suffice to post the resulting BVP and provide the details essential for deriving the stabilized finite-element method. Furthermore, we limit the discussion to homogeneous Dirichlet boundary conditions for notational simplicity. Consequently, the BVP can be summarized as follows.
Let $\Omega$ be an open, bounded, and polyhedral domain of $\mathds{R}^2$ and $\partial \Omega$ the domain's boundary. We want to find $u$ such that
\begin{subequations}
\begin{align}\label{eq:reynolds}
    -\tfrac{1}{12}\nabla\cdotn \big(H^3(x,y)\,\nabla \left(g(u)\,\uk(x,y)\right)\big)
        -\partial_x\big(\left(g(u)-1\right)H(x,y)\,\uk(x,y)\big)
        &=f(x,y) -\partial_x H(x,y) \quad && \text{in}\ \Omega, \\
    u &= 0 \quad &&\text{on\ } \partial\Omega.
\end{align}
\end{subequations}
The \textit{gap function} $H(x,y)$ is assumed to be known, typically as a result of a previous simulation step. Hence, the problem can be considered quasi-static,\footnote{In general, terms involving explicit time-dependency may also appear in the Reynolds equation~\eqref{eq:reynolds} in some applications, but we will not discuss the transient cases in this work.} i.e., we aim to find a solution $\uk(x,y)$ for a given $H(x,y)$. In our numerical examples, we will assume the gap function to be of the form 
\begin{equation}\label{eq:gapFunction}
    H  =  1-\zeta\,\cos(x-x_a)
\end{equation}
 with constants $\zeta$ and $x_a$. Furthermore, $g(u)$ is a 
\textit{switch function}, describing the transition between the pressure and cavitation zones. Here, we will use the regularized version~\cite{Nitzschke2016a}
\begin{equation}
    g(\uk)=\frac{1}{\uppi} \arctan \left(\frac{\uk}{1-\overline{\uk}}\right)+\frac{1}{2},
\end{equation}
where the constant $\overline{\uk}$ is usually chosen in the range $0.9< \overline{\uk}<1$.
For other regularization approaches in cavitation algorithms, see, e.g.,~\cite{Fesanghary2011,bayada2001,Zhang2017,nilsson2007}.
A few remarks may help to clarify the BVP in this particular setting:
\begin{itemize}\setlength\itemsep{-0.2\baselineskip}
    \item Physically, the computational domain represents a thin cylindrical layer (assuming a journal bearing) and is often formulated in cylindrical coordinates $(z,\theta)$. As the radial dependency and the curvature are commonly neglected, the computational domain reduces to a plane surface. Thus, we write the PDE in Cartesian coordinates to avoid any confusion regarding the differential operators. 
    \item The above version of the Reynolds equation is already non-dimensionalized; see \cite{Pfeil2023} for details.
    \item The continuous unknown field $\uk(x,y)$ has different physical interpretation according to its sign. In regions where $\uk \geq 0$ (`pressure zone'), it represents the hydrodynamic pressure, while it is related to a film fraction\footnote{The film fraction (also known as density ratio)~$\vartheta$ is a measure of the relative amount of liquid in the cavitated fluid film.} where $\uk < 0$ (`cavitation zone'). More specifically, the non-dimensionalized pressure $p$ and the film fraction $\vartheta$ can be derived from $\uk$ as $p=g\uk$ and $\vartheta=(1-g)\uk+1$, respectively.\footnote{The idea of describing the complementary unknowns~$p$ and~$\vartheta$ by a common global unknown~$u$ was proposed by Shi and Paranjpe~\cite{Shi2002}. In their work, they use~$p$ and~$\vartheta$ in the strong form of the differential equation and substitute these quantities by~$u$ after discretization, meaning that~$p$ and~$\vartheta$ are both interpolated separately, and the substitution with~$u$ is performed afterwards with respect to the discrete nodal quantities. In contrast, we interpret~$u$ as a continuous field constituting the global unknown in our strong form. Alternatively, it is possible to employ concepts of mixed finite element formulations and treat $p$ and~$\vartheta$ as separate unknown fields \cite{Lengiewicz2014}. }
    \item The first term of the PDE is dominant in the pressure zone; the second is dominant in the cavitation zone. The regularization function ensures a smooth transition between both regimes.
    \item The function $f(x,y)$ incorporates any additional forcing terms. In this study, we solely included it for the purpose of creating manufactured solutions.
    \item The signs in Eq.~\eqref{eq:reynolds} are chosen such that the diffusion term in the weak form will be positive definite (In contrast to, e.g., \cite{Pfeil2023}), as this is a common assumption in the analysis of similar problems.
\end{itemize}
Introducing shorthand notations for the left-hand side and right-hand side of \eqref{eq:reynolds}, the PDE is of the generic form 
\begin{equation}\label{eq:strongFormGeneric}
    \mathcal{L}(u,u) = \hat{f}.
\end{equation}
Note that $\mathcal{L}(u,u)$ is a quasi-linear operator; its second slot accounts for the spatial derivatives and the first one for the dependence of their coefficients on the unknowns.

\section{Variational form}\label{variational}\noindent
We denote by $\mathcal{W}\subset H^1_0(\Omega)$ the adequate function space where the continuous problem is well posed, which depends on the choice of $g(u)$. Here, the standard notation is used for the Sobolev spaces, i.e., $H^m(\Omega)$ is the space of functions whose distributional derivatives up to integer order $m \geq 0$ are square-integrable in $\Omega$. The space of square-integrable functions is denoted as $L^2(\Omega)$, and $H_0^1(\Omega)$ consists of functions in $H^1(\Omega)$ that vanish on $\partial \Omega$. 

Multiplying Eq.~\eqref{eq:reynolds} by test functions $\tf \in \mathcal{W}$ (i.e., taken from the same space as $\uk$) and integrating over the computational domain yields
\begin{equation}\label{eq:weak}
    -\tfrac{1}{12} 
    \big( \tf,\, \nabla\cdotn(H^3\,\nabla (g\,\uk)) \big)
    - \big(\tf,\, \partial_x (\left(g\shortm1\right)H\uk )\big)
    = \big(\tf,\, f- \partial_x H \big).
\end{equation}
Here and in the following, we use the symbol $(\dotP,\dotP)$ to indicate integration over the computational domain, irrespective of whether or not this expression represents the $\mathrm{L}^2$-inner product in $\Omega$.
Integrating the first term by parts, we obtain
\begin{equation}\label{eq:weak2}
    \tfrac{1}{12}\big(\nabla\tf,\, H^3\nabla (g\uk)\big)
        -\big(\tf,\, \partial_x ((g\shortm1 )H\uk )\big)
        =\big(\tf,\, f- \partial_x H \big).
\end{equation} 
It is convenient to define a function $\phi(\uk)=g(\uk)\,\uk$, such that
\begin{equation}\label{eq:definePhi}
    \nabla (g(\uk)\,\uk) =  \frac{\partial}{\partial \uk}(g(\uk)\,\uk)\,\nabla\uk=\phi'(\uk)\,\nabla\uk,
\end{equation}
where the prime symbol denotes a derivative with respect to $\uk$. Substituting Eq.~\eqref{eq:definePhi} into the first term of Eq.~\eqref{eq:weak2} and expanding the second term leads to
\begin{equation}\label{eq:weak3}
    \tfrac{1}{12}\big(\nabla\tf,\, H^3\phi'\,\nabla\uk\big)
        -\big(\tf,\, (g\shortm1)H\, \partial_x \uk\big)
        -\big(\tf,\, \partial_x (\left(g\shortm1\right)H)\uk\big)
        =\big(\tf,\, f-\partial_x H \big).
\end{equation}
The above weak form resembles that of the standard convection-diffusion-reaction (CDR) equation, and we will use this terminology throughout the paper, regardless of the physical meaning of these terms. However, observe that all terms on the left-hand side involve nonlinearities (since $g=g(u)$) and spatial variation of the coefficients ($H=H(x,y)$). We may also note that, in this particular setting, the `convection' term involves only derivatives with respect to $x$.

For easier comparison with previous developments for the (linear) CDR equation, let us introduce the abbreviations 
\begin{align*}
    &k(u,x,y) = \tfrac{1}{12}H^3(x,y)\,\phi'(u),\quad 
    &&\bm{a}(u,x,y) = [(g(u)-1)\,H(x,y),\,0]^\T,\\ 
    &s(u,x,y) = \partial_x \big((g(u)\shortm1)\,H(x,y)\big), \quad 
    &&\hat{f}(x,y) = f(x,y)-\partial_x H(x,y). 
\end{align*}
Hence, the weak form of the problem reads: Find $u\in \mathcal{W}$ such that
\begin{equation}\label{eq:weakAbb}
    \big(\nabla\tf,\, k(u,x,y) \nabla\uk\big)
        -\big(\tf, \bm{a}(u,x,y) \cdotn \nabla \uk\big)
        -\big(\tf,\, s(u,x,y)\, \uk\big)
        =\big(\tf, \hat{f}(x,y) \big)
\end{equation}
for all $v\in \mathcal{W}$. More compactly, we may write
\begin{equation}\label{eq:bilinearFormEquation}
    B(u; u, v) = L(v),
\end{equation}
where, for a given $\mhat{u}\in \mathcal{W}$, we denote by $B(\mhat{u};\dotP, \dotP)$ the bilinear form defined on $\mathcal{W}\times \mathcal{W}$ as
\begin{equation}\label{eq:B}
    B(\mhat{u}; u, v) = \big(\nabla\tf,\, k(\mhat{u},x,y) \nabla\uk\big)
        -\big(\tf,\, \bm{a}(\mhat{u},x,y) \cdotn \nabla \uk\big)
        -\big(\tf,\, s(\mhat{u},x,y)\, \uk\big), 
\end{equation}
and $ L(v)$ is given as 
\begin{equation}\label{eq:L}
    L(v) = \big(\tf, \hat{f}(x,y) \big).
\end{equation}

\begin{rmk}\label{rmk:artificialDiffusivity}
    As mentioned in the introduction, the most popular remedy to avoid unphysical oscillations in the solution of the Reynolds equation consists in including a term representing artificial diffusion (or using an upwind scheme that adds this diffusion indirectly), thus ensuring that the problem remains well-posed even in the cavitation region where $k(u,x,y)$ (nearly) vanishes. For the solution to converge to the correct result under mesh refinement, the artificial diffusivity must decrease with the element size $h$. The value of the artificial diffusivity can be justified, e.g., by the streamline diffusion method~\cite{Lengiewicz2014}, which is an upwind scheme based on a Petrov-Galerkin formulation. This approach leads to the form
    \begin{equation}\label{eq:artificialdiffusion}
        D(\mhat{u}; u, v) =  \big(\hat{\bm{a}}\cdot\nabla\tf,\, \tfrac{h}{2} \nabla\cdot (\bm{a}\,\uk)\big)
        = -\big(\partial_x\tf,\, \tfrac{h}{2} \partial_x (a_x\,\uk)\big)
    \end{equation}
    with the convection vector $\bm{a}=[a_x,\,0]^\T$ defined above and $\hat{\bm{a}} = \bm{a}/|\bm{a}|$.
    Note that, in the case of the Reynolds equation, the artificial diffusivity takes significant values only in the cavitation zone. As shown in~\cite{Pfeil2023}, the term derived in Eq.~\eqref{eq:artificialdiffusion} is also exactly equivalent to the numerical diffusion obtained by the upwind difference technique, which is commonly used for the stabilization of FDM and FVM solutions of the Reynolds equation. Here, we will only use this artificial diffusion formulation for comparison and validation of our results in the numerical examples.
\end{rmk}

\section{Stabilized finite element formulation}\label{stabilization}\noindent
In this work, we consider only two-dimensional plane geometries and hence require a polygonal finite element partition of the computational domain $\Omega$. We may assume the partition to be quasi-uniform with an element diameter denoted as $h$. 
We construct the finite element spaces $\mathcal{W}_h \subset \mathcal{W}$ in the usual manner employing polynomial interpolation of uniform order throughout the computational domain. Assuming the same finite element spaces for the trial and test functions, we obtain the Galerkin discretization as
\begin{equation}\label{eq:weakFEM}
    B(u_h; u_h, v_h) = L(v_h).
\end{equation}
It is well-known that the discretized CDR equation is prone to nonphysical oscillations in the case of convection-dominated problems. The same issue is observed when solving the Reynolds equation. In fact, the influence of cavitation within the model employed here always leads to the problem being `convection'-dominated in a cavitation domain. Note that, without regularization, the diffusion term would exactly vanish in the cavitation domain. 

We choose a stabilization analogously to the linear CDR-equation that has been studied extensively in previous work \cite{Codina1998,Codina2000a,Codina2000,Principe2010}. A particularly detailed discussion of the approach employed here is presented in \cite{Principe2008}. The underlying idea is based on the variational multiscale (VMS) method, which relies on splitting the unknown function $u$ into the component that belongs to the chosen finite element space and a remainder (referred to as subscale and denoted by $\til{u}$) that cannot be resolved by the finite element discretization. The corresponding space of subscales is denoted as $\til{\mathcal{W}}$. Thus, we have 
\begin{equation}
    u = u_h + \til{u}, \quad \mathcal{W} = \mathcal{W}_h \oplus \til{\mathcal{W}}.
\end{equation}
Substituting the above decomposition into the variational form and separating equations with respect to the test functions yields
\begin{subequations}\label{eq:stabilizedWFseparated}
    \begin{alignat}{6}
        \label{eq:stabilizedWFseparatedGT}
        &B(\ukh ;\ukh, \tfh) &\ +\  &B(\ukh ;\til{\uk}, \tfh)\, &= L(\tfh), \\ 
        \label{eq:stabilizedWFseparatedSS}
        &B(\ukh ;\ukh, \til{\tf}) &\ +\ &B(\ukh ;\til{\uk}, \til{\tf})\, &= L(\til{\tf})\,.
    \end{alignat}
    \end{subequations}
Note that we used $B(\ukh;\dotP, \dotP)$, hereby computing the coefficients of the weak form based on the finite element approximation $\ukh$. An alternative, referred to as \textit{nonlinear subscales}, consists in including the contributions from the subscales in the computation of the coefficients, i.e., using instead $B(u;\dotP, \dotP)$. 
Following the standard procedure of the VMS method (see, e.g., \cite{Codina2017a}), the second term in Eq.~\eqref{eq:stabilizedWFseparatedGT} is integrated by parts. In this step, it is a common assumption that the subscales vanish on the element boundaries. Thus, we obtain
\begin{equation}\label{eq:stabilizedWFadjoint}
    B(\ukh ;\ukh, \tfh) + \sum_K \big( \til{\uk}, \mathcal{L}^*(\uk,\tfh)\big) _K \approx L(\tfh),
\end{equation}
where the subscript $K$ indicates integration and summation over all elements in the mesh. In the above equation, $\mathcal{L}^*$ is the adjoint of the differential operator $\mathcal{L}$. Several different approaches exist within this family of VMS methods, which mainly differ in the approximation of the subscale equation \eqref{eq:stabilizedWFseparatedSS}.
As an important subclass of methods, residual-based approaches assume a subscale model of the form
\begin{equation}\label{eq:subscaleEvolution}
    \tau^{\shortminus 1}\,\til{\uk} = \til{P}\left( \hat{f} - \mathcal{L}(\uk,\ukh) \right)
\end{equation}
with a stabilization parameter $\tau$ that approximates the effect of the differential operator on the subscales. Furthermore, $\til{P}$ denotes the projection onto the subscale space \cite{Castillo2019}, here applied to the residual of the finite-element approximation. In this work, we follow the concept of \textit{orthogonal subgrid scales} (OSGS) \cite{Codina2000}, i.e., we choose $\til{P}=P^\perp$ to be the projection onto the space orthogonal to the finite element space. It can be rewritten as $P^\perp=I-P_h$, where $P_h$ is the projection onto the finite element space. 
Substituting Eq.~\eqref{eq:subscaleEvolution} into \eqref{eq:stabilizedWFadjoint} then yields
\begin{equation}\label{eq:stabilizedSubsti}
    B(\ukh ;\ukh, \tfh) + \sum_K \big( \tau\,P^\perp( \hat{f} - \mathcal{L}(\uk,\ukh)), \mathcal{L}^*(\uk,\tfh)\big) _K \approx L(\tfh).
\end{equation}
However, when employing orthogonal subgrid scales, we can drastically simplify the above formulation by considering only those terms in the residual as well as the adjoint operator $\mathcal{L}^*$ that are essential for achieving stability. 
As discussed in detail in the literature on the OSGS approach, the orthogonal projection of the relevant terms is sufficient for achieving stability while maintaining optimal convergence~\cite{Codina2008, Castillo2019}. In our application, where we want to stabilize convection, it suffices to include the convection term in $\mathcal{L}$ and $\mathcal{L}^*$. Hence, we can expect to obtain a stable formulation by approximating the stabilization terms as\footnote{Keep in mind that $\bm{a} = \bm{a}(u,x,y)$ and $\tau = \tau(u,x,y)$, which is omitted in the following equations for readability. Hence, the stabilization term is nonlinear. }
\begin{equation}\label{eq:stabi}
    \sum_K \big( \tau\,P^\perp\negthickspace\left( \hat{f} - \mathcal{L}(\uk,\ukh) \right), \mathcal{L}^*(\uk,\tfh)\big) _K 
    \approx\big(\bm{a}\cdotn\nabla\tfh,\, \tau\, P^\perp(\bm{a}\cdotn \nabla \ukh)\big) \eqqcolon S(\ukh;\, \ukh, \tfh)
\end{equation}
or, considering $P^\perp=I-P_h$,
\begin{equation}\label{eq:stabTerm}
    S(\ukh;\,\ukh, \tfh) = 
    \big(\bm{a}\cdotn\nabla\tfh,\, \tau\, \bm{a}\cdotn \nabla \ukh\big) 
    - \big(\bm{a}\cdotn\nabla\tfh,\, \tau\, P_h\big(\bm{a}\cdotn \nabla \ukh\big)\big).
\end{equation}
The first term of the stabilization $S(\ukh;\,\tfh,\ukh)$ can be interpreted as additional nonlinear diffusion and leads to a stable method. The second term, involving the projection onto the finite-element space, ensures optimal convergence under mesh refinement. Details on the computation of this projection are provided in the next section.

The stabilization parameter $\tau$ is adapted from previous work on the linear CDR equation \cite{Codina2000,Codina2000a}. It needs to be chosen such that the residual of the governing equation vanishes as the element size $h$ approaches zero, and the optimal convergence rate is attained in the asymptotic regime. A common approach to deriving a computable expression for the stabilization parameter is based on a Fourier analysis of the residual, see, e.g., \cite{Codina2002,Principe2010,Castillo2019}. From the discussion in the papers mentioned before, it is known that the stabilization parameter for the linear CDR equation, say  $-\hat{k}\,\Delta u + \hat{\*a}\cdot\nabla u + \hat{s}\,u=0$, can be chosen of the form $(c_1\,\hat{k}/h^2 + c_2\,|\hat{\*a}|/h + \hat{s})^{\shortminus1}$. Here, we need to take into account the nonlinear and inhomogeneous nature of $k(u,x,y)$, $\bm{a}(u,x,y)$, and $s(u,x,y)$; thus, $\tau(u,x,y)$ will be evaluated at every Gauss point when integrating the stabilization term \eqref{eq:stabi}. We obtain
\begin{align} \nonumber
    \tau(u,x,y)
    &=\left(\frac{c_1}{h^2} \left|k(u,x,y)\right| + \frac{c_2}{h} \left| \bm{a}(u,x,y)\right| + \left|s(u,x,y)\right|\right)^{-1} \\
    &= \left(\frac{c_1}{12\,h^2} \left|H^3(x,y)  \phi'(\uk)\right| + \frac{c_2}{h} \left| (g(u)\shortm1)H(x,y)\right| + \left|\partial_x \big((g(u)\shortm1)\,H(x,y)\big)\right|\right)^{-1}, \label{eq:tau} 
\end{align}
where $c_1,\,c_2$ are algorithmic constants commonly chosen (for linear elements) as 
\begin{equation}
    c_1=4,\quad c_2=2.
\end{equation}
In summary, the stabilized finite element method reads:  Find $\ukh\in \mathcal{W}_h$ such that
\begin{equation}\label{eq:stabFEMfinal}
    B(\ukh ;\ukh, \tfh) + S(\ukh;\,\tfh,\ukh)  =  L(\tfh),
\end{equation}
for all $\tfh\in \mathcal{W}_h$ with $B$, $L$, and $S$ given in Eqs.~\eqref{eq:B}, \eqref{eq:L}, and \eqref{eq:stabTerm} and the stabilization parameter $\tau$ defined in \eqref{eq:tau}.

\section{Computing the projection}\label{projection}\noindent
The projections of the convection term present in the stabilized weak form are apriori unknown functions in the finite-element space; we will denote them as 
\begin{equation}
    \xi_h = P_h\big(\bm{a}\cdotn \nabla \ukh\big).
\end{equation}
Computing these functions $\xi_h$ corresponds to solving
\begin{equation}\label{eq:projection}
    \big(\eta_h,\, \bm{a}\cdotn \nabla \ukh\big) - \big(\eta_h, \xi_h \big) = 0
\end{equation}
with test functions $\eta_h\in \mathcal{W}_h $ (i.e., they consist of the same basis functions as $\tfh$).
Combining Eqs.~\eqref{eq:stabFEMfinal} and \eqref{eq:projection}, we obtain the system of equations
\begin{subequations}\label{eq:system}
    \begin{align}
    B_s(\ukh ;\ukh, \tfh) 
    - \big(\bm{a}\cdotn \nabla \tfh,\, \tau\, \xi_h \big)
    &=L(\tfh)\label{eq:system_a}\\ 
    \big(\eta_h,\, \bm{a}\cdotn \nabla \ukh\big) - \big(\eta_h, \xi_h \big) &= 0, \label{eq:system_b}
\end{align}
\end{subequations}
where we introduced, for notational convenience,
\begin{equation}
    B_s(\ukh; \tfh,\ukh) = 
    B(\ukh;\ukh, \tfh) 
    + \big(\bm{a}\cdotn\nabla\tfh,\, \tau\, \bm{a}\cdotn \nabla \ukh\big).
\end{equation}
In this work, we compute the projections $\xi_h$ \textit{implicitly} along the lines of the approach suggested in \cite{Codina2008}.
That is to say, at a given iteration $i$, the solution $\ukh^{i+1}$  as well as the projections $\xi_h^{i+1}$ are both treated as unknowns and obtained by solving the following coupled system of equations:
\begin{subequations}\label{eq:system_imp}
    \begin{align}
        B_s(\ukh^i;\tfh,\ukh^{i+1})
    - \Big(\bm{a}\cdotn \nabla \tfh,\, \tau\, \xi_h^{i+1} \Big)
    &=L(\tfh), \label{eq:system_imp_a}\\ 
    \Big(\eta_h,\, \bm{a}\cdotn \nabla \ukh^{i+1} \Big) 
    - \Big(\eta_h, \xi_h^{i+1} \Big) 
    &= 0. \label{eq:system_imp_b}
\end{align}
\end{subequations}
Details on the linearization of the Galerkin terms in $B_s(\ukh^i;\tfh,\ukh^{i+1})$ will be presented in Section~\ref{linearization}.
For clarity, we rewrite the system \eqref{eq:system_imp} in matrix form, assuming, for instance, a finite-element version using standard piece-wise polynomial shape functions:
\begin{equation}
    \left[\begin{array}{ll}
        \m{K} & -\m{P}_\tau \\
        \m{P} & -\m{M}
    \end{array}\right] \left[\begin{array}{rr}
        \vt{u}\\
        \vg{\xi} 
    \end{array}\right]
    = \left[\begin{array}{rr}
        \vt{F}\\
        \vt{0} 
    \end{array}\right].
\end{equation}
Here, the symbols $ \vt{u}$ and $\vg{\xi}$ denote the vectors of unknowns, i.e., the coefficients of the finite-element representation (typically nodal values) of $\ukh^{i+1}$ and $\xi_h^{i+1}$. The definition of the matrices $\m{K},\,\m{P}_\tau,\,\m{P},\,\m{M}$ follows directly from comparing with Eqs.~\eqref{eq:system_imp}.\footnote{In some applications of orthogonal subgrid scales, the stabilization parameter $\tau$ is assumed constant in the entire domain (hence $\m{P}_\tau = \tau\, \m{P}^\T $) or approximated as being constant within each element. In our case, $\tau = \tau(u,x,y)$ is a rather complicated function of the solution $u$ as well as the spatial coordinates; hence, it is evaluated at every Gauss point.} 
Thus, in the implicit version of the linearization, we formally introduced additional degrees of freedom representing the projection of the stabilization term onto the finite-element space. Nevertheless, these additional degrees of freedom can be eliminated by static condensation, yielding
\begin{equation}
    \left(\m{K} - \m{P}_\tau \m{M}^{\shortm1}\m{P}\right)\vt{u} = \vt{F}\,.
\end{equation}
When using a finite element version that allows \textit{mass lumping} -- i.e., the Gram matrix $\m{M}$ can be diagonalized -- the cost of computing the additional term $\m{P}_\tau \m{M}^{\shortm 1}\m{P}$ is small compared to the solution of the final system of equations. 

We note that alternatives to this implicit approach are generally possible. In particular, a straightforward explicit version consists in replacing ${i+1}$ by $i$ in Eq.~\eqref{eq:system_imp_b} and the second term in Eq.~\eqref{eq:system_imp_a}, thus using the projections at the previous iteration $i$ when computing $\uk^{i+1}$. This approach is essentially a staggered scheme involving alternating computations of $\uk$ and $\xi$. However, it is known from numerical studies on the current and previous applications that this explicit scheme leads to poor convergence of the nonlinear terms. In transient problems, on the other hand, the usually preferred option relies on using the projections at the previous \textit{time step}; hence, they are not updated by the nonlinear solver and, consequently, do not affect the convergence of the same. The resulting error in the projections is typically negligible for sufficiently small time steps. Such an approach may also prove beneficial for the solution of the Reynolds equation when it is part of a larger system in which the gap function $H$ varies in time. However, for the current work where we are interested in the stationary case, we will stick to the implicit version outlined above. 

\section{Linearization of the Galerkin terms}\label{linearization}\noindent
As all terms in the weak form \eqref{eq:bilinearFormEquation} are nonlinear, we must choose a suitable linearization and employ an iterative solver for obtaining a converged solution. The simplest approach is obtained by a fixed-point iteration (also referred to as Picard's method):
\begin{equation}
    B_{\mathrm{P}}(u^{i}; u^{i+1}, v) = 
    \big(\nabla\tf,\, k(u^{i},x,y) \nabla\uk^{i+1}\big)
    -\big(\tf,\, \bm{a}(u^{i},x,y) \cdotn \nabla \uk^{i+1}\big)
    -\big(\tf,\, s(u^{i},x,y)\, \uk^{i+1}\big).
\end{equation}
This approach simply corresponds to using the result $\uk^i$ of a previous iteration $i$ in the integration of the stiffness matrix for computing the solution $\uk^{i+1}$. 
Analogously, the stabilization terms are linearized by approximating
\begin{multline} 
    S_{\mathrm{P}}(u^i;\,\ukh^{i+1},\tfh) = 
    \big(\bm{a}(u^i,x,y)\cdotn\nabla\tfh,\, \tau(u^i,x,y)\, \bm{a}(u^i,x,y)\cdotn \nabla \ukh^{i+1}\big) \\
    - \big(\bm{a}(u^i,x,y)\cdotn\nabla\tfh,\, \tau(u^i,x,y)\, P_h\big(\bm{a}(u^i,x,y)\cdotn \nabla \ukh^{i+1}\big)\big). \label{eq:stabTermPicard}
\end{multline}
In addition, we make use of the classical Newton-Raphson method in order to obtain significantly faster convergence when the initial guess is sufficiently close to the converged solution. To this end, we approximate the diffusive term by a first-order Taylor expansion at $\uk^i $:
\begin{align*}
    k(u)\,\nabla\uk 
    &\approx k(\uk^i)\,\nabla\uk^i 
    + \left.\partial_\uk\left(k(\uk)\,\nabla\uk\right)\right|_{\uk^i} (\uk^{i+1}-\uk^i)
    \\
    & = k(\uk^i)\,\nabla\uk^i 
    + k'(\uk^i)\,\nabla\uk^i\,(\uk^{i+1}-\uk^i) + k(\uk^i)\,\nabla (\uk^{i+1}-\uk^i)
    \\
    & = k(\uk^i)\,\nabla \uk^{i+1}
    +k'(\uk^i)\,\nabla\uk^i\,\uk^{i+1} 
    -k'(\uk^i)\,\nabla\uk^i\,\uk^i.
\end{align*}
For convenience, we rewrite the convective and reactive terms (the second term in Eq.~\eqref{eq:weak}) as
\begin{equation}
    \partial_x \big(\left(g\shortm1\right)H\uk\big) 
    = \partial_x \big(\psi(\uk)\,H\big) 
    %=  \psi(\uk)\,\partial_x H + H\,\partial_x \psi(\uk)
    =  \psi(\uk)\,\partial_x H + H\, \psi'(\uk)\,\partial_x\uk
\end{equation}
with $\psi(\uk) = (g-1)\,\uk$, and we obtain an analogous linearization as 
\begin{align*}
    &\psi(\uk)\,\partial_x H + H\, \psi'(\uk)\,\partial_x\uk\\
    &\approx
    \psi(\uk^i)\,\partial_x H + H\, \psi'(\uk^i)\,\partial_x\uk^i 
    + \left.\partial_\uk \left( \psi(\uk)\,\partial_x H + H\, \psi'(\uk)\,\partial_x\uk\right)\right|_{\uk^i} (\uk^{i+1}-\uk^i)\\
    &=
    \psi(\uk^i)\,\partial_x H + H\, \psi'(\uk^i)\,\partial_x\uk^i 
    + \left( 
    \psi'(\uk^i)\,\partial_x H 
    + H\, \psi''(\uk^i)\,\partial_x\uk^i 
    \right)(\uk^{i+1}-\uk^i)
    + H\, \psi'(\uk^i)\,\partial_x 
     (\uk^{i+1}-\uk^i)\\
     &=
     H\, \psi'(\uk^i)\,\partial_x 
     \uk^{i+1}  
    + \left( 
    \psi'(\uk^i)\,\partial_x H 
    + H\, \psi''(\uk^i)\,\partial_x\uk^i 
    \right)(\uk^{i+1}-\uk^i)
    + \psi(\uk^i)\,\partial_x H.
\end{align*}
Hence, the linearization of the (unstabilized) weak form \eqref{eq:bilinearFormEquation} based on the Newton-Raphson method reads
\begin{multline}
    B_\mathrm{N}^{i+1}(\tf,\,\uk) = \left(\nabla\tf,\, k^{\,i}\,\nabla \uk^{i+1}\right)
    +\left(\nabla\tf,\, k'^{\,i}\,\nabla\uk^i\,\uk^{i+1} 
    \right)
    -\left(\nabla\tf,\, k'^{\,i}\,\nabla\uk^i\,\uk^i\right)
    \\
    -\left(\tf,\, H\, \psi'^{\,i}\,\partial_x \uk^{i+1}\right)
   - \left(\tf,\, \chi^i\uk^{i+1} \right)
   +\left(\tf,\, \chi^i\uk^i \right)
   - \left(\tf,\,  \psi^{\,i}\,\partial_x H\right) = L(v)
\end{multline}
with the abbreviations
\begin{equation*}
    \chi^i \coloneqq 
    \psi'^{\,i}\,\partial_x H(x,y) 
    + H(x,y)\, \psi''^{\,i}\,\partial_x\uk^i, \qquad
    k^i \coloneqq k(\uk^i,x,y), \qquad
    \psi^i \coloneqq \psi(\uk^i).
\end{equation*}
Implementing higher-order approximations of the stabilization terms $S(\tfh,\ukh)$ and the involved projections is cumbersome, mainly due to the relatively complicated dependency of the stabilization parameter $\tau$ on the solution $\uk$. We will adhere to using the linearization by Picard's method for the stabilization terms as described above (Eq.~\eqref{eq:stabTermPicard}) even when employing Newton's method for the standard Galerkin terms. This simplification is common in the vast majority of works on subgrid-scale methods, and our numerical results confirm the rapid convergence of this approximated Newton-Raphson scheme.

\section{Shock-capturing}\label{shockCapturing}\noindent
The stabilization outlined before yields a stable method that converges at an optimal rate and avoids unphysical global oscillations. On the other hand, local oscillations may still arise in confined regions where the solution changes abruptly. While this behavior is typical in boundary layers, it can, in our current application, also occur in the transition between cavitation and pressure zones, where the behavior of the solution changes significantly. Such oscillations manifest themselves in `over- and undershooting' phenomena, which are ubiquitous in many interpolation problems. The presence of such oscillations is not necessarily problematic, and their magnitude decays under mesh refinement in the stabilized method. Nevertheless, it can sometimes be desired to suppress oscillations and obtain a smoother solution, particularly on relatively coarse meshes.
The corresponding techniques are often referred to as `shock-capturing,' and their further advancement is still an active research area not restricted to stabilized finite-element methods. Essentially, many of these techniques are based on adding local artificial diffusion of a magnitude that depends on the gradient of the solution. Hence, the shock-capturing term is usually nonlinear, and its magnitude must be carefully selected in order to yield smoother results without deviating too far from the correct solution \cite{Castillo2014, Badia2014}.
Here, we will not delve into the details of different shock-capturing approaches but simply apply one established version with proven capabilities to reduce local oscillations without significantly deteriorating accuracy. The idea -- as introduced in \cite{Codina1993} and slightly modified in \cite{Knopp2002} -- 
is to introduce artificial diffusion that is scaled by the strong residual of the governing PDE. Hence, we add a term of the form
\begin{equation}
    D(\hat{u}_h;\ukh,\tfh)= \sum_K \big(\nabla\tfh,\, \tau_s(\hat{u}_h) \nabla \ukh\big)_K 
\end{equation}
with 
\begin{equation}
    \tau_s(u) = R_K^*(u)  \, \sigma_K(u)
\end{equation}
and 
\begin{equation}
    R_K^*(u) = \frac{|| \mathcal{L}(u,u) - \hat{f}||_{L^2(K)}}{\alpha\,||\hat{f}||_{L^2(K)} + ||u||_{H^1(K)}}, \qquad
    \sigma_K(u) =  \frac{h_K}{2}\,\operatorname{max}\negthinspace\left(0,\,\beta-\frac{1}{P_K}\right), \qquad
    P_K = \frac{h_K\,R_K^*(u)}{2\,k(u)},
\end{equation}
where the subscript $K$ again indicates evaluation on each element. Hence, the artificial diffusivity is scaled by the element-wise $L^2$ norm of the strong residual, adequately normalized by the element-wise $H^1$ norm of the solution and the element-wise $L^2$ norm of the right-hand side. The factor $\alpha$ is introduced here for dimensional consistency and chosen as 
\begin{equation}
    \alpha = \frac{\sum_K ||u||_{H^1(K)}}{\sum_K ||\hat{f}||_{L^2(K)}}.
\end{equation} 
Furthermore, $\sigma_K(u)$ is a limiter function, ensuring positive bounded values and yields small artificial diffusivity in regions where the existing diffusivity $k(u)$ is large. The algorithmic constant $\beta$ is chosen as 0.7 in our numerical examples. For details on this approach, the reader is referred to~\cite{Codina1993}.

\section{Numerical examples}\label{numex}

\subsection{Convergence test}\label{sec:convergence}
\begin{figure}[tbh]\centering
    \subfloat[\label{fig:manu_u}]{\includegraphics[width=0.50\textwidth]{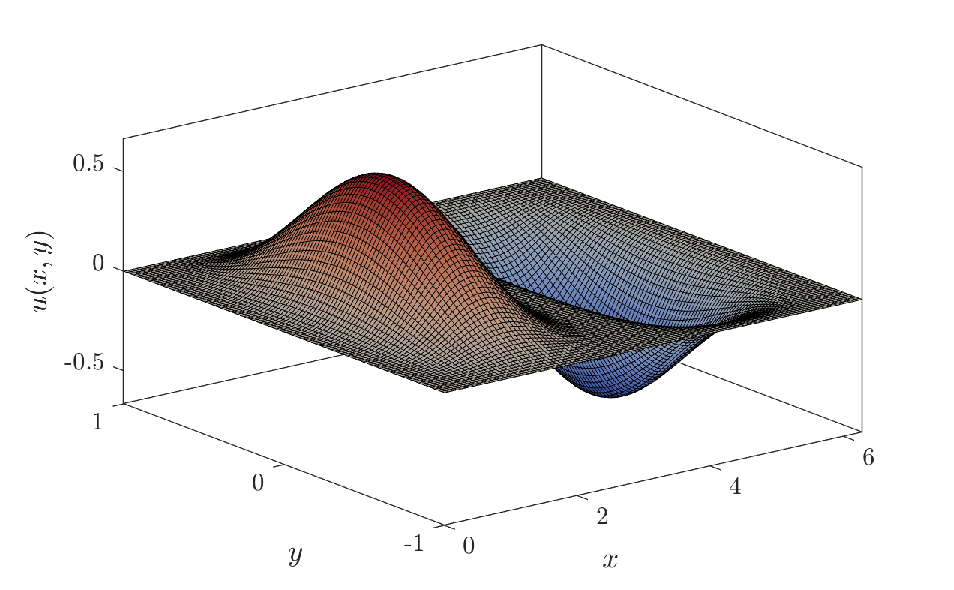}}\hfill
    \subfloat[\label{fig:manu_f}]{\includegraphics[width=0.50\textwidth]{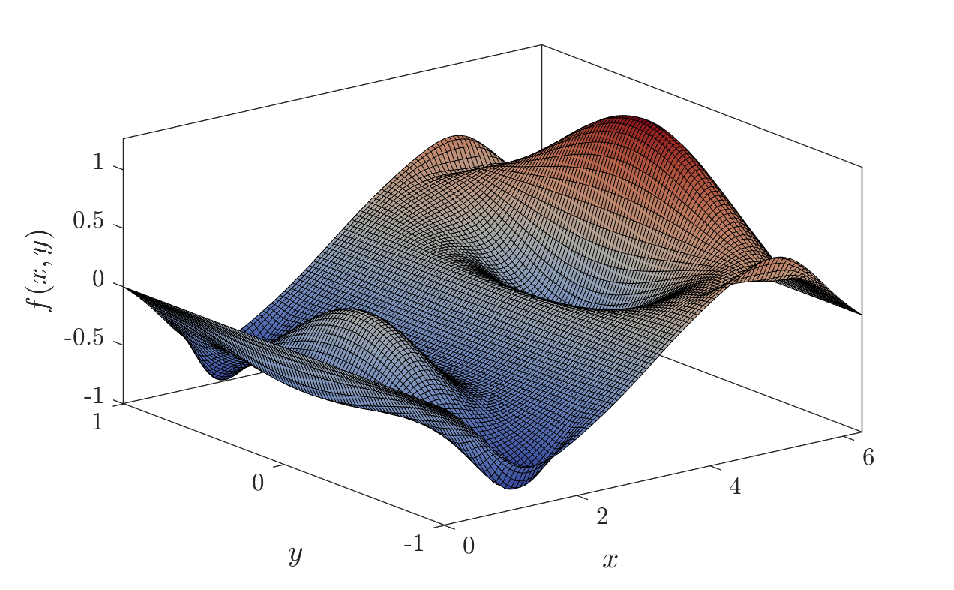}}
    \caption{Manufactured solution (a), forcing function (b) for the example analyzed in Section~\ref{sec:convergence}.}
\end{figure}

\begin{figure}\centering
    \subfloat[with stabilization \label{fig:manu_uVSx_stab}]{\includegraphics[width=0.5\textwidth]{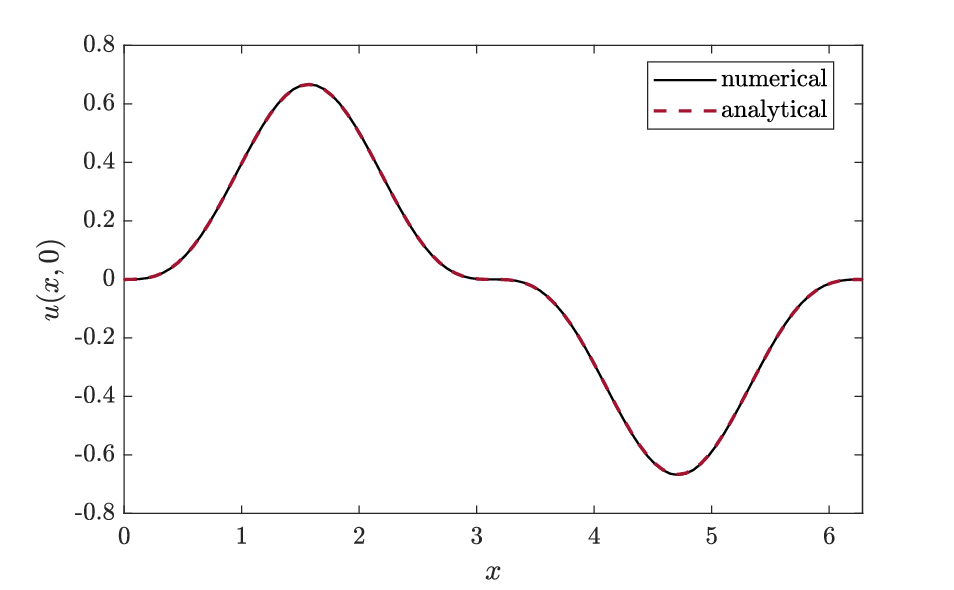}}
    \subfloat[without stabilization \label{fig:manu_uVSx_nostab} ]{\includegraphics[width=0.5\textwidth]{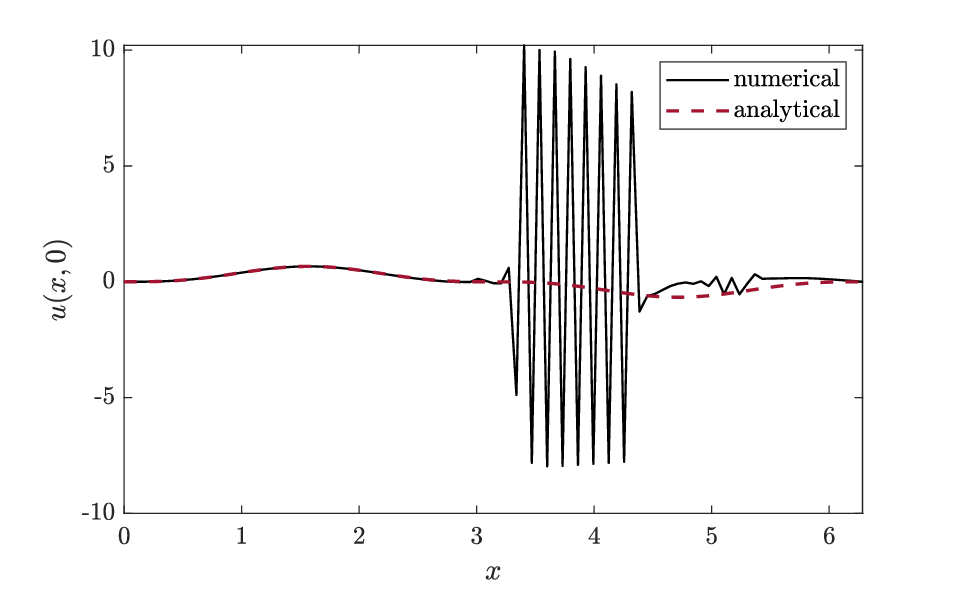}}
    \caption{Analytical and numerical solution (for the finest mesh of $96\times 32$ elements) along the line $y=0$. \label{fig:manu_uVSx} }
\end{figure}

\begin{figure}[tbh]\centering
    \subfloat[\label{fig:manu_taua}]{\includegraphics[width=0.50\textwidth]{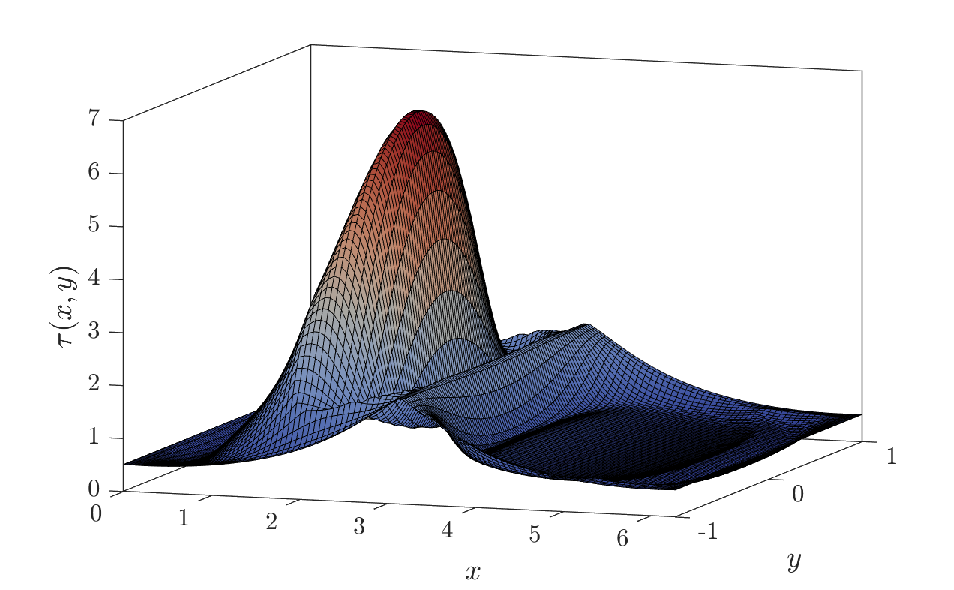}}\hfill
    \subfloat[\label{fig:manu_taub}]{\includegraphics[width=0.50\textwidth]{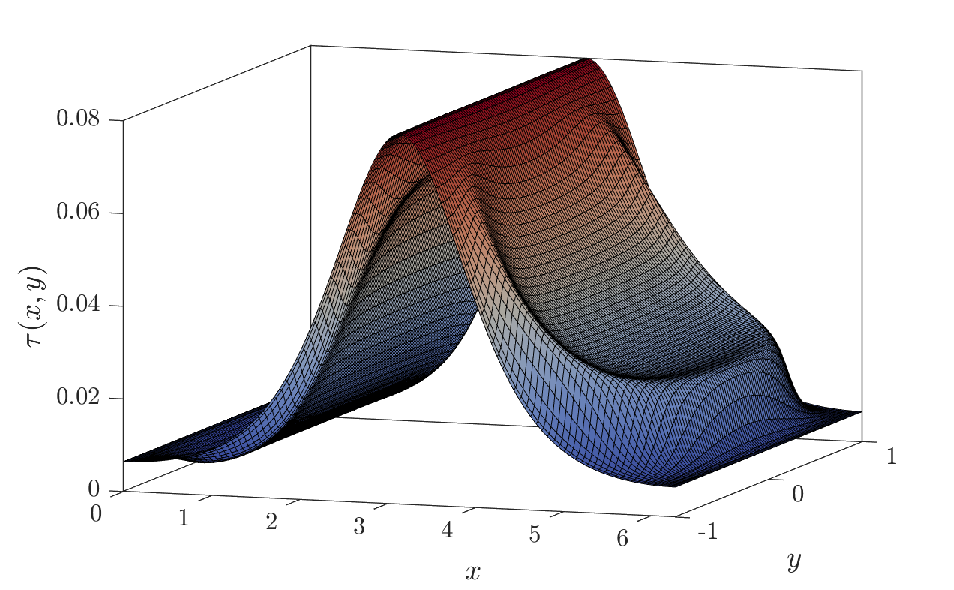}}
    \caption{Stabilization parameter $\tau(x,y)$ computed for the coarsest (a) and finest (b) mesh in the example \ref{sec:convergence}. \label{fig:manu_tau}}
\end{figure}

\noindent
We begin the numerical studies by assessing a simple example with a known smooth solution in the computational domain
\begin{equation*}
    \Upomega=\left\{\left. (x,\,y) \in \mathbb{R}^2\ \right|\ 0\leq x\leq 2\uppi,\ \shortminus 1\leq y\leq 1\right\}.    
\end{equation*}
We construct this problem using a so-called manufactured solution. That is to say, we choose $u(x,y)$ and substitute this function into the governing PDE to obtain the corresponding forcing term. In our example, we use
\begin{equation}\label{eq:manuCosSin}
    u(x,y) =  (1-\cos(2x)) \sin(x)(1+\cos(\uppi\,y))/6
\end{equation} 
as depicted in Fig.~\ref{fig:manu_u}. Note that the solution takes positive and negative values (i.e., involves a pressure zone and a cavitation zone) and vanishes at the domain's boundary. 
The gap function is chosen according to Eq.~\eqref{eq:gapFunction} with the parameters $\zeta = \tfrac{1}{2}$, $x_a = \uppi$. In order to ensure the method's robustness, different values of the regularization parameter $\overline{\uk}$ have been tested, namely $\overline{\uk} = 0.9,\,0.91,\,...\,,\,0.99$, leading to almost identical results. The solutions presented in this section have been computed using $\overline{\uk}=0.98$. 
Substituting Eq.~\eqref{eq:manuCosSin} into the governing PDE yields the forcing function $f(x,y)$ depicted in Fig.~\ref{fig:manu_f}. We refrain from printing the lengthy expressions describing this function; however, a Matlab function for computing $f(x,y)$ for this example is publicly available \cite{Gravenkamp2023e}. Homogeneous Dirichlet conditions are prescribed at the boundary, and the initial guess is chosen as $u_0(x,y) = 1$ inside the domain.
The computational domain is discretized by starting with three quadrilateral elements along the $x$-direction (of size $\sfrac{2}{3}\,\uppi \times 2$) and consecutively dividing each element into four to obtain a series of consistently refined meshes.
As an example, Fig.~\ref{fig:manu_uVSx_stab} displays the computed solution along the line $y=0$ when using a mesh of $96\times 32$ elements. The comparison against the analytical solution shows no significant discrepancies when employing the stabilized finite element method. In contrast, if we do not include the stabilization term (but leave all other parameters unchanged), we are unable to find a converged solution. Figure~\ref{fig:manu_uVSx_nostab} shows the solution after 50 fixed-point iterations, revealing strong oscillations in the cavitation region, which prevent convergence in this example. However, it should be noted that this behavior is strongly mesh-dependent, and in this simple example with a smooth solution, the fixed-point iteration did converge for some of the meshes we tested. To get a notion of the amount of stabilization included in the method, Fig.~\ref{fig:manu_tau} shows the variation of the stabilization parameter $\tau(x,y)$ within the entire domain for the coarsest and finest mesh utilized in this example.

\begin{figure}
    \subfloat[error \label{fig:manu_convergence}]{\includegraphics[width=0.5\textwidth]{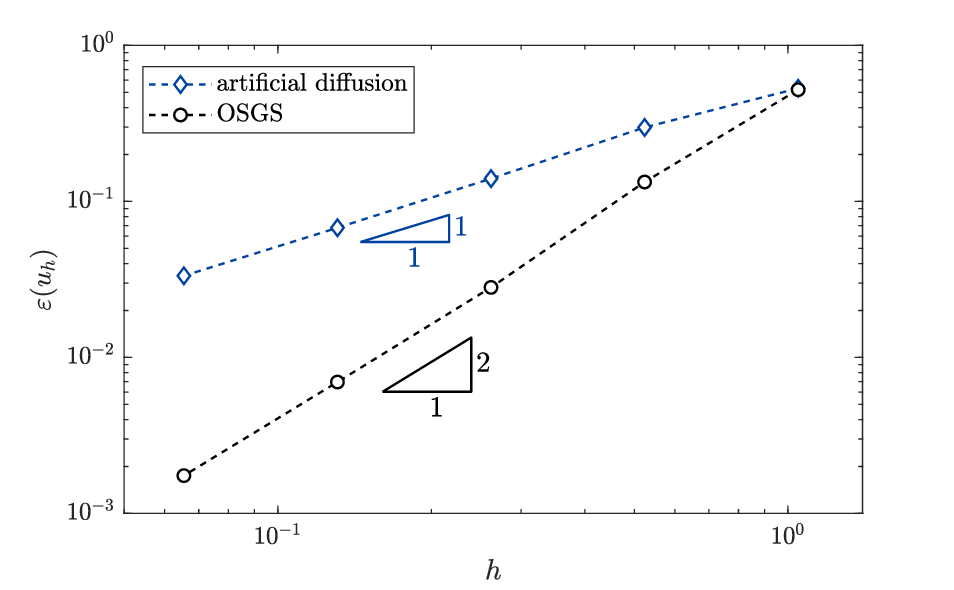}}\hfill
    \subfloat[residual \label{fig:manu_residual}]{\includegraphics[width=0.5\textwidth]{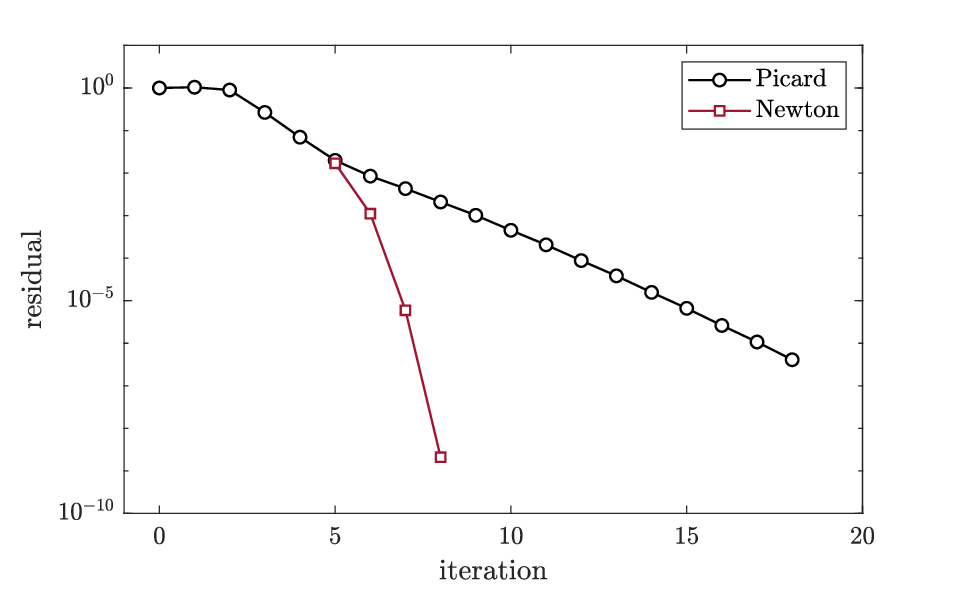}}
    \caption{Convergence of the error when performing $h$-refinement (a); convergence of the residual when using Picard's iteration or Newton's method (b). }
\end{figure}

To better assess the accuracy of the stabilized method, we compute the normalized $L^2$ norm of the error for different mesh sizes, i.e.,
\begin{equation}
    \epsilon(u_h) = \frac{||u-u_h||_{L^2}}{||u||_{L^2}}.
\end{equation}
This error measure decreases proportionally to $h^2$, which is the optimal convergence rate in the case of linear elements; see Fig.~\ref{fig:manu_convergence}. For comparison, we show the corresponding results when employing artificial diffusion (as described in Remark~\ref{rmk:artificialDiffusivity}) instead of the proposed OSGS approach. While the artificial diffusion does allow the method to converge, it leads to only first-order convergence and overall significantly larger errors. However, it should be acknowledged that the OSGS approach causes slightly larger computational costs for the same mesh, as it requires the computation of the additional matrices related to the projections (cf.~Eq.~\eqref{eq:system_imp}), as well as their static condensation.

In Fig.~\ref{fig:manu_residual}, we depict the convergence of the residual when employing fixed-point iteration and Newton's method in the proposed approach -- again for the discretization using $96\times 32$ elements. We obtain the expected asymptotic convergence rates of approximately one and two, respectively, when averaging the last three points of the depicted graphs. Note that Newton's method tends to fail when starting with an initial guess far from the true solution. For this reason, we usually perform a few fixed-point iterations (in this example four) before switching to Newton's method for faster convergence.

\subsection{Example involving large gradients}\label{sec:boundaryLayer}
\begin{figure}
    \subfloat[\label{fig:boundaryLayer_u}]{\includegraphics[width=0.50\textwidth]{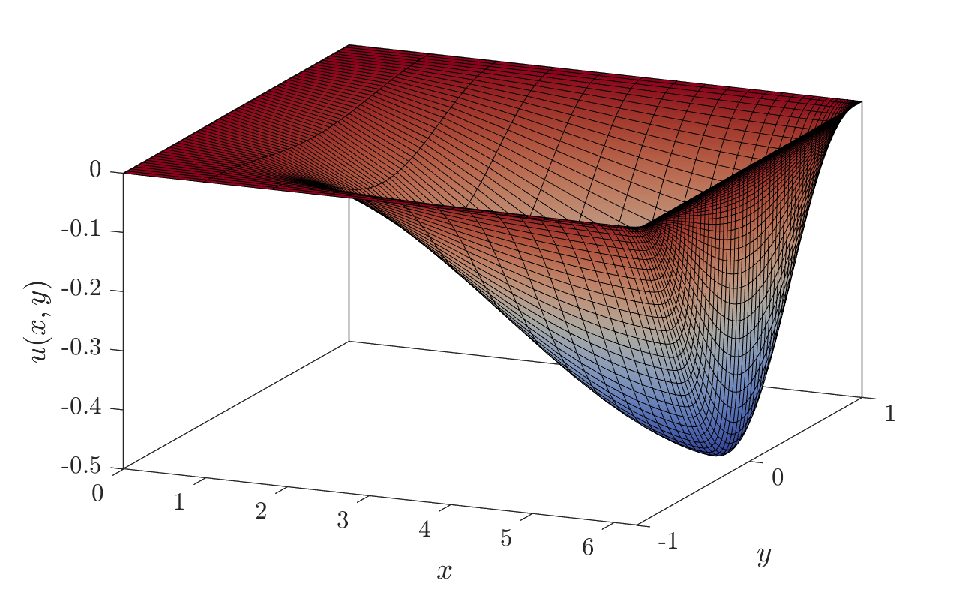}}\hfill
    \subfloat[\label{fig:boundaryLayer_f}]{\includegraphics[width=0.50\textwidth]{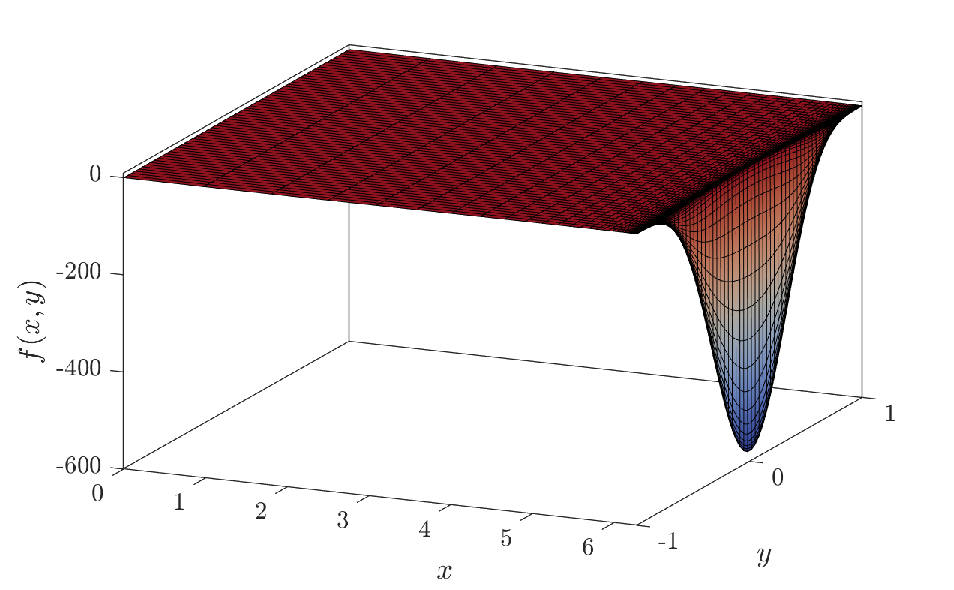}}
    \caption{\label{fig:boundaryLayer}Manufactured solution (a) and forcing function (b) for the example analyzed in Section~\ref{sec:boundaryLayer}.}
\end{figure}

\begin{figure}\centering
    \subfloat[with stabilization \label{fig:boundaryLayer_uVSx_stab}]{\includegraphics[width=0.5\textwidth]{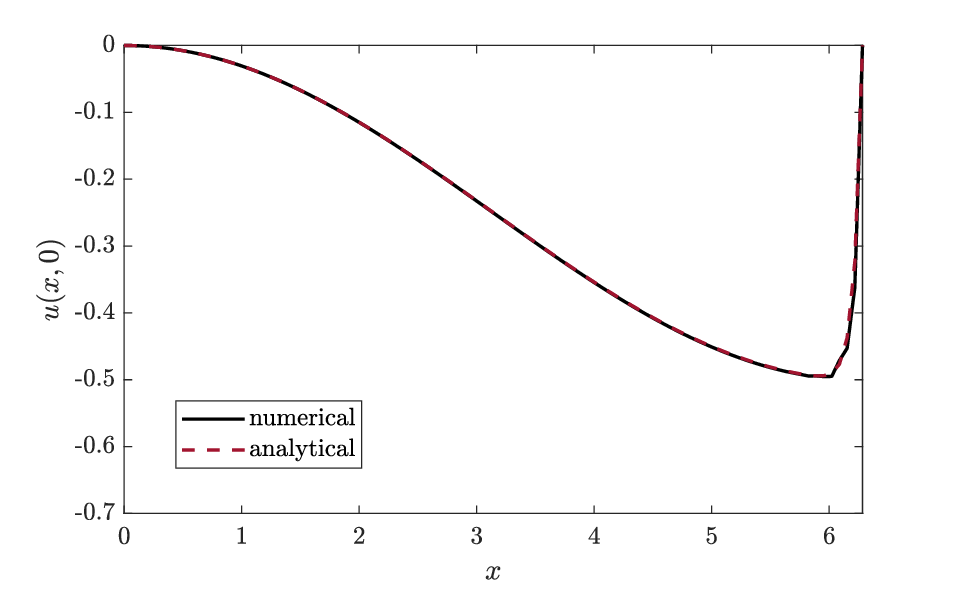}}\hfill
    \subfloat[without stabilization \label{fig:boundaryLayer_uVSx_nostab} ]{\includegraphics[width=0.5\textwidth]{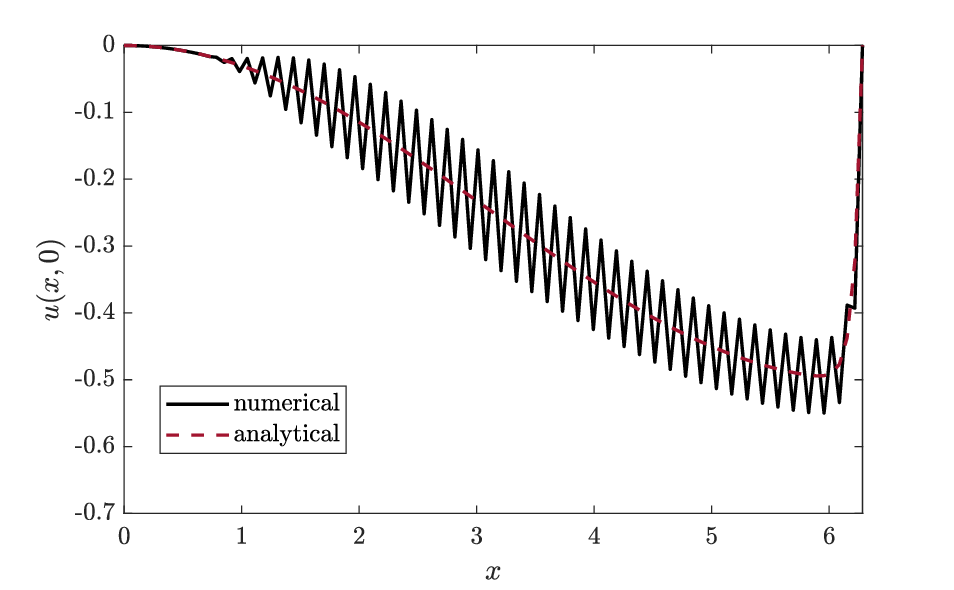}}
    \caption{Analytical and numerical solution for the finest mesh of $96\times 32$ elements along the line $y=0$. \label{fig:boundaryLayer_uVSx} }
\end{figure}

\begin{figure}\centering
    \subfloat[with shock capturing \label{fig:boundaryLayer_uVSx_nosc}]{\includegraphics[width=0.5\textwidth]{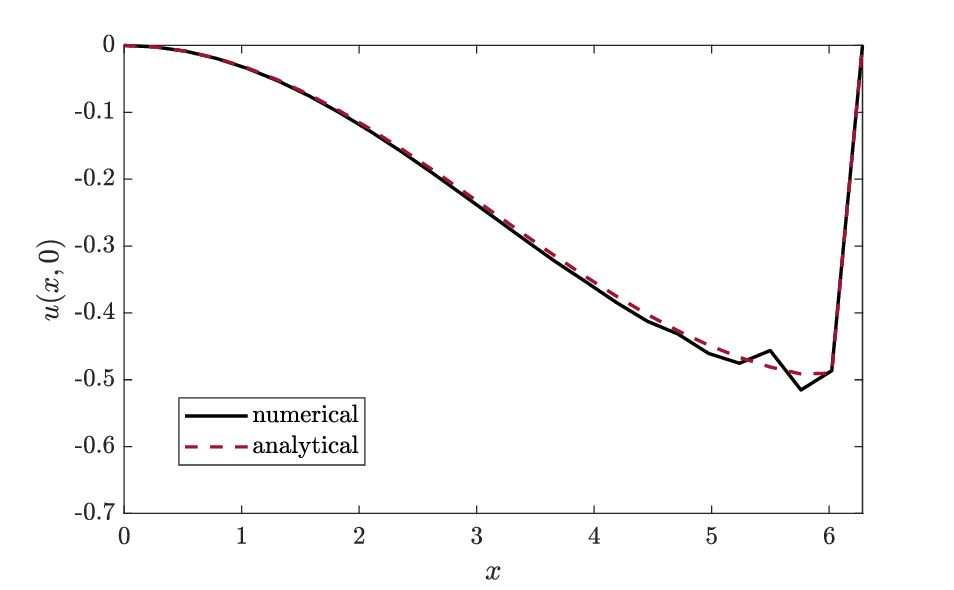}} \hfill
    \subfloat[without shock capturing \label{fig:boundaryLayer_uVSx_sc}]{\includegraphics[width=0.5\textwidth]{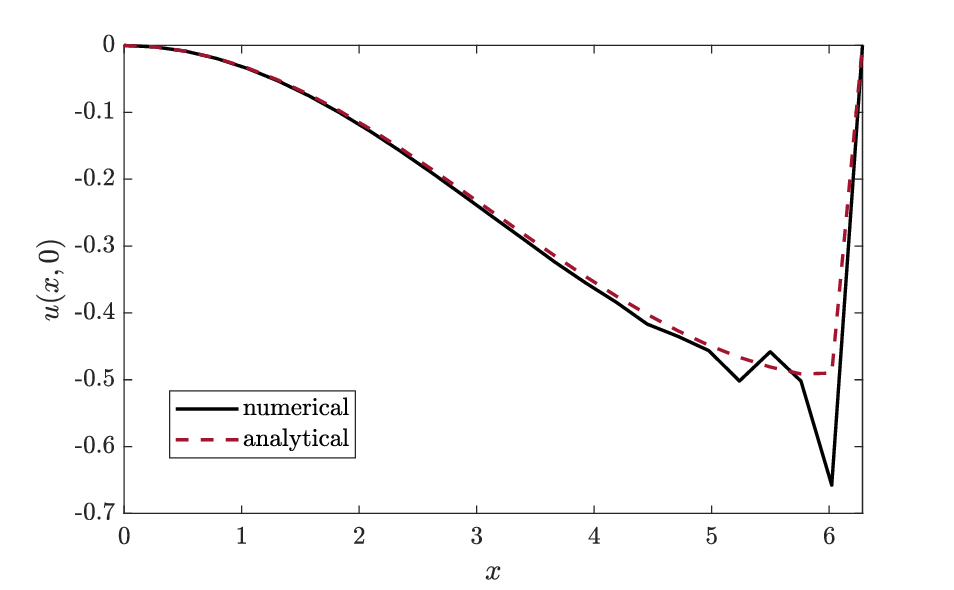}}
    \caption{Analytical and numerical solution for a coarse mesh of $24\times 8$ elements along the line $y=0$. \label{fig:boundaryLayer_uVSx_h3} }
\end{figure}

\begin{figure}
    \subfloat[error \label{fig:boundaryLayer_convergence}]{\includegraphics[width=0.5\textwidth]{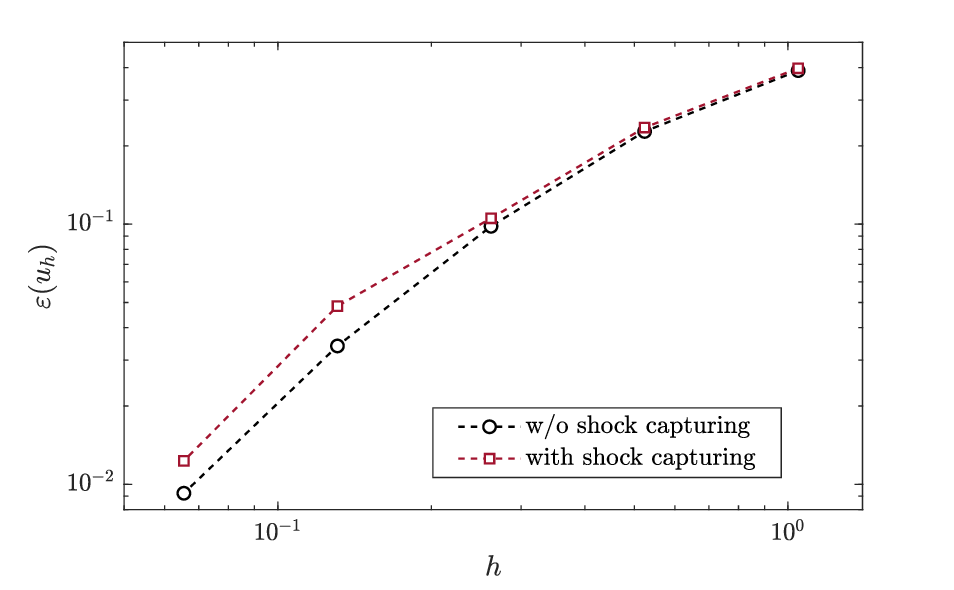}}\hfill
    \subfloat[residual \label{fig:boundaryLayer_residual}]{\includegraphics[width=0.5\textwidth]{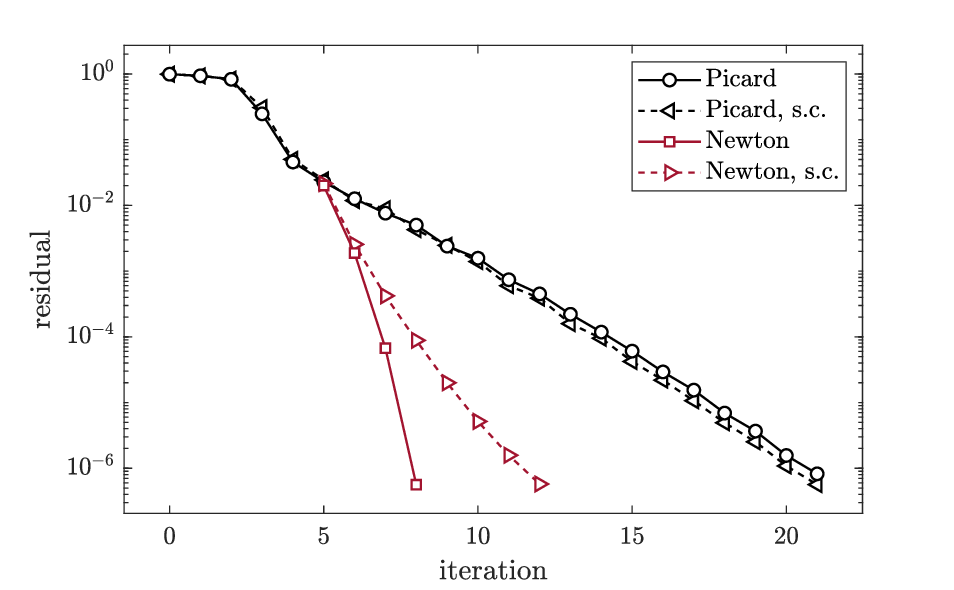}}
    \caption{Convergence of the error when performing $h$-refinement (a); convergence of the residual when using Picard's iteration or Newton's method (b). Results are presented with and without shock-capturing (s.c.).}
\end{figure}

\noindent
As a second example, we consider another manufactured solution that involves large gradients both in the field $u(x,y)$ as well as the forcing term $f(x,y)$, see Fig.~\ref{fig:boundaryLayer}. 
The solution is chosen as 
\begin{equation}
    u(x,y) = \frac{1}{4}\Big(\frac{1-\operatorname{e}^{\frac{100}{2\uppi} x}}{1-\operatorname{e}^{100}}-1+\tfrac{1}{2}(\cos(\tfrac{x}{2})+1)\Big) (1+\cos(\uppi y)).
\end{equation}
The computational domain, gap function, regularization, as well as the boundary conditions, initial guess, and meshes, are identical to the previous example in Section~\ref{sec:convergence}. Again, modifying the regularization parameter has little effect on the convergence of the method, and the presented results are computed using  $\overline{\uk} = 0.98$.
The solution along the line $y=0$ when using the mesh of $96\times 32$ elements is depicted in Fig.~\ref{fig:boundaryLayer_uVSx}. When employing the stabilized method, this discretization yields a sufficiently accurate solution. In contrast, severe oscillations are still visible without stabilization. 
The effect of incorporating the shock-capturing term is demonstrated using a coarser mesh of $24\times 8$ elements in Fig.~\ref{fig:boundaryLayer_uVSx_h3}. 
While the stabilization leads to optimally convergent results and avoids global oscillations, local oscillations of large amplitude can still occur near the boundary where large gradients are present in the exact solution. In this case, the shock-capturing approach helps reduce these oscillations. On the other hand, note that the stabilized method gives accurate results on sufficiently fine meshes even without shock-capturing -- at least in this simple example. The graph in Fig.~\ref{fig:boundaryLayer_uVSx} was computed without shock-capturing, and there is no visible difference when incorporating it. It should also be noted that shock-capturing can have undesired effects on the convergence both with respect to mesh refinement and convergence of the nonlinear terms. In the approach we used, this effect depends strongly on the parameter $\beta$. Large values of $\beta$ lead to smoother results but increase the overall error of the numerical solution and may deteriorate the convergence of the nonlinear solver. Here, we chose a value of $\beta = 0.7$, which has been found to yield a suitable trade-off between effective shock-capturing and accuracy. From the results of the convergence study in Fig.~\ref{fig:boundaryLayer_convergence}, we can observe that the shock-capturing then leads to only a slight increase in the overall error. Regarding the convergence of the nonlinear solver, we mainly see an effect on Newton's method; see Fig.~\ref{fig:boundaryLayer_residual}. This can be attributed to the fact that the nonlinear coefficient (involving element-wise integration of the strong residual) must be computed based on the previous iteration; hence, this term is still linearized by a simple fixed-point iteration, even when using Newton's method for all other terms. The results in this figure are again evaluated for the finest mesh, but those for coarser meshes show a similar trend.

\subsection{Realistic example} \label{sec:realExample}
\begin{figure}\centering
    \subfloat[with shock-capturing \label{fig:real_SC}]{\includegraphics[width=0.50\textwidth]{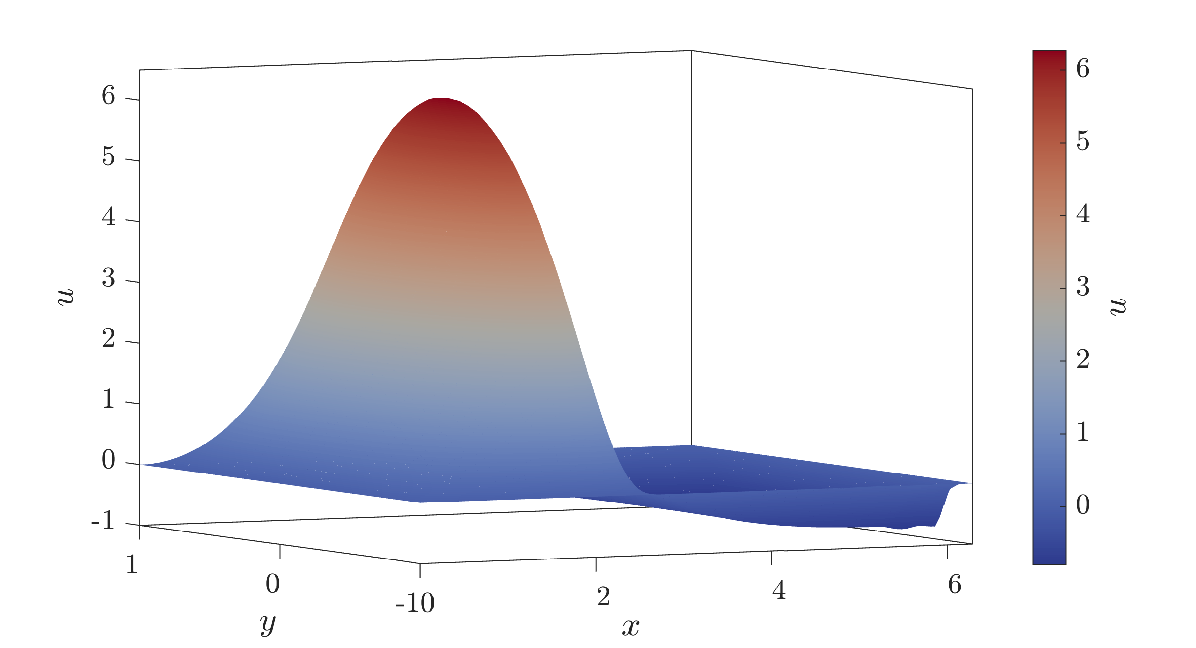}}\hfill
    \subfloat[without shock-capturing \label{fig:real_noSC}]{\includegraphics[width=0.50\textwidth]{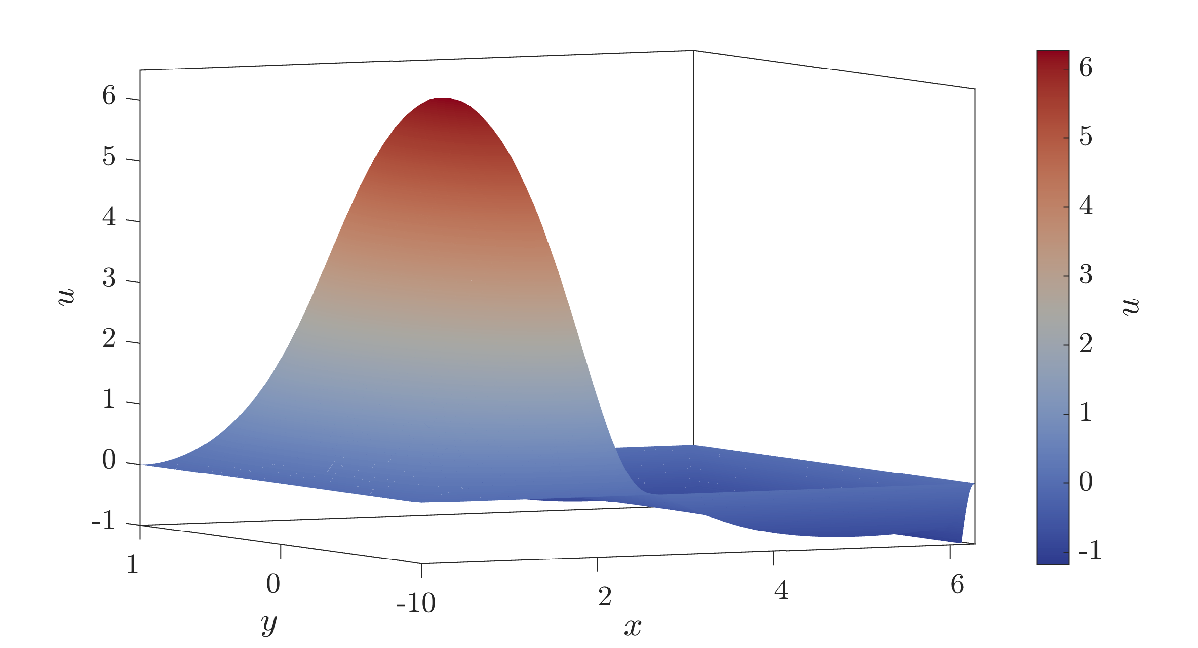}}
    \caption{Solution of the example analyzed in Section~\ref{sec:realExample} with and without shock-capturing. \label{fig:real_u}}
\end{figure}

\begin{figure}\centering
    \subfloat[\label{fig:real_uVSx_all}]{\includegraphics[width=0.5\textwidth]{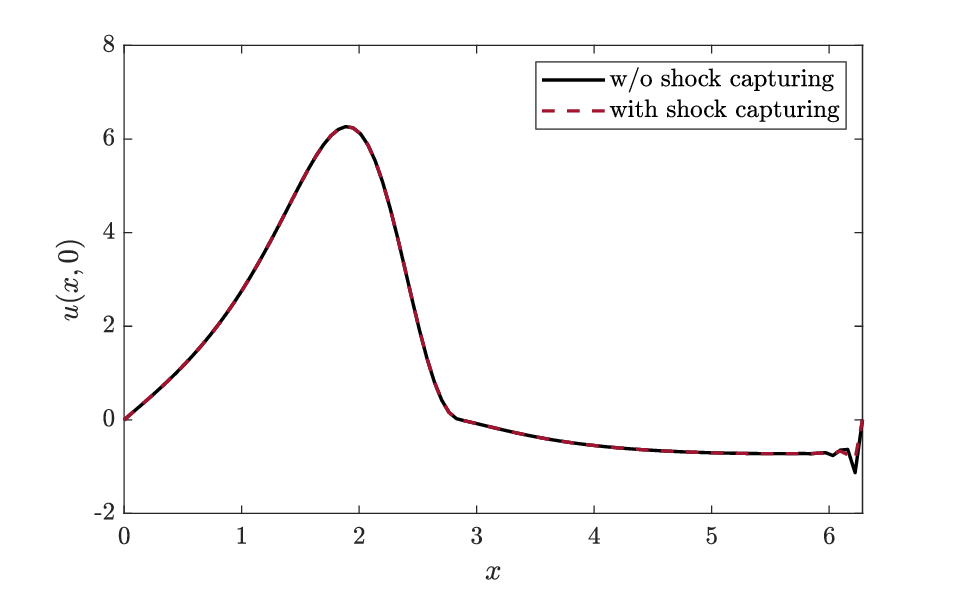}}\hfill
    \subfloat[\label{fig:real_uVSx_zoom}]{\includegraphics[width=0.5\textwidth]{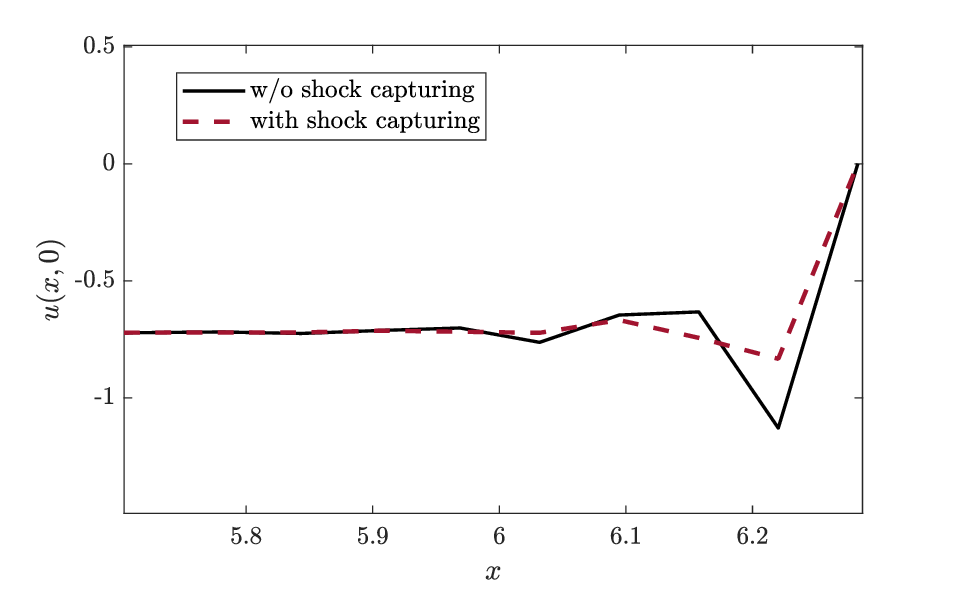}}
    \caption{Solution with and without shock-capturing; (a) along the line $y=0$, (b) detail near the boundary.\label{fig:real_uVSx}}
\end{figure}

\begin{figure}\centering
    \includegraphics[width=0.5\textwidth]{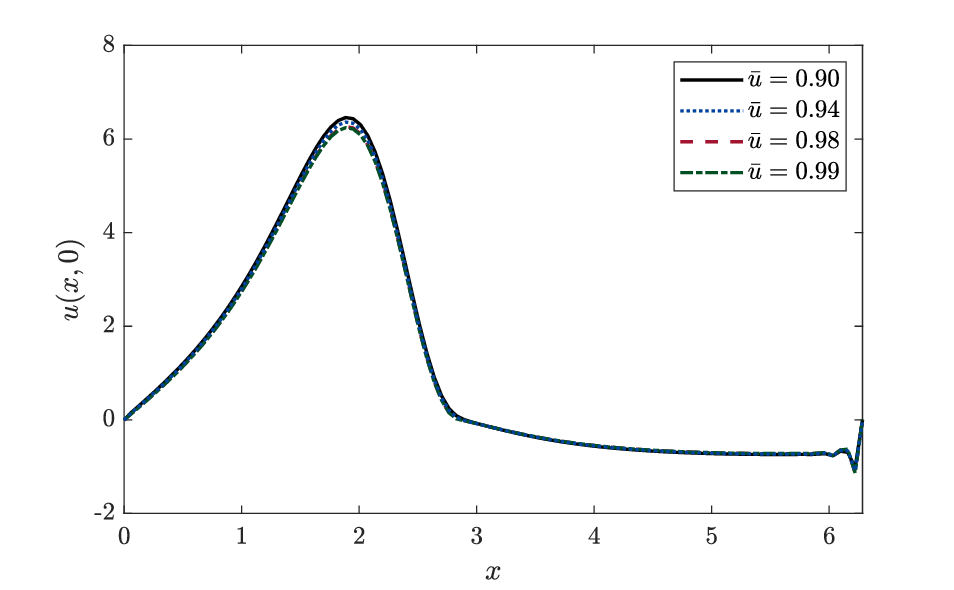}
    \caption{Solution along the line $y=0$ for different values of the regularization parameter $\overline{\uk}$.\label{fig:real_uVSx_ubar}}
\end{figure}

\noindent
As a final numerical study, we consider a more realistic scenario of a bearing, leading to a typical distribution of pressure and cavitation zones. The example is based on one described in detail in \cite{Pfeil2023}. Only the boundary conditions are simplified here for conciseness.  
The computational domain is the same as in the previous examples (i.e., a rectangle of dimensions $2\uppi\times 2$). The gap function is defined by the parameters $\zeta = 0.6$, 
$x_a = \tfrac{7}{9}\uppi$ ($140^\circ$ in a cylindrical coordinate system). Homogeneous Dirichlet boundary conditions are applied. The initial guess is again chosen as $u_0(x,y)=1$. Figure~\ref{fig:real_u} shows the solution in the entire computational domain, computed using a mesh of $100\times 32$ elements. The result displays the typical smooth behavior inside the pressure zone and a rather abrupt transition to the cavitation zone. For a more in-depth discussion of the underlying physics and the relevance within the scope of bearing simulations, we refer the interested reader to \cite{Pfeil2023}. Similar to the example in Section~\ref{sec:boundaryLayer}, the solution features large gradients near the boundary $x=2\uppi$. This gives rise to oscillations that are mitigated when employing the shock-capturing approach. Figure~\ref{fig:real_uVSx} shows the solution along the line $y=0$. It can be observed that the proposed stabilized method avoids global oscillations inside the domain. The shock-capturing approach successfully suppresses the local oscillations near the boundary without causing any visible deviation elsewhere in the computational domain. These results were once again computed using a regularization parameter of $\overline{\uk}=0.98$. In addition, we present in Fig.~\ref{fig:real_uVSx_ubar} the solution along $y=0$ using different values of $\overline{\uk}$. Only small deviations can be noticed for rather low values of  $\overline{\uk}\leq 0.94$, confirming that the model is robust for a sufficient range of reasonable values.

\section{Conclusion}\noindent
The proposed stabilized finite element method for the Reynolds equation follows rather straightforwardly from the existing approaches available for the CDR equation when taking into account the nonlinearities inherent in the cavitation model. We showed that it is sufficient to include one additional term in order to stabilize convection and demonstrated numerically the optimal convergence of the resulting method. Furthermore, we have found that local oscillations, which can occur in regions exhibiting large gradients, can effectively be suppressed by existing shock-capturing techniques.

\section*{Acknowledgements}\noindent
H.~Gravenkamp acknowledges grant CEX2018-000797-S funded by the Ministerio de Ciencia e Innovaci\'on, MCIN/AEI/10.13039/501100011033. S.~Pfeil acknowledges the financial support by the German Research Foundation, Germany (DFG), project no.\ 490625563, grant no.\ \mbox{WO~2085/8-1}.
R.~Codina acknowledges the support received from the ICREA Acad\`emia Research Program of the Catalan Government.

\bibliographystyle{elsarticle-num}
\bibliography{reynoldsOSS.bib,reynoldsOSS_Pfeil.bib}

\end{document}